\documentclass[letterpaper, 10 pt, conference]{ieeeconf}  % Comment this line out if you need a4paper

\IEEEoverridecommandlockouts                              % This command is only needed if 
\usepackage{cite}
\usepackage{amsmath,amssymb,amsfonts}
\usepackage{graphicx}
\usepackage{textcomp}
\usepackage{xcolor}
\usepackage[english]{babel}
\usepackage{tabularx}
\usepackage{algorithm}
\usepackage{algpseudocode}
\usepackage{subcaption}

\newtheorem{assumption}{Assumption}
\newtheorem{definition}{Definition}

\def\BibTeX{{\rm B\kern-.05em{\sc i\kern-.025em b}\kern-.08em
    T\kern-.1667em\lower.7ex\hbox{E}\kern-.125emX}}

\title{\LARGE \bf
Collision-Free Continuum Deformation Coordination of a Multi-Quadcopter System Using Cooperative Localization
}

\author{Hamid Emadi$^{1}$, Harshvardhan Uppaluru$^{1}$, Hashem Ashrafiuon$^{2}$, and  Hossein Rastgoftar$^{1}$
\thanks{$^{1}$H. Emadi, H. Uppaluru, and H. Rastgoftar are  with the Aerospace and Mechanical Engineering Department at University of Arizona, AZ, USA, Emails: \{hamidemadi, huppaluru, hrastgoftar\}@arizona.edu}
\thanks{$^{2}$H. Ashrafiuon is  with the Mechanical Engineering Department at Villanova University, PA, USA 19085, Email: hashem.ashrafiuon@villanova.edu}
}

\begin{document}

\maketitle
\thispagestyle{empty}
\pagestyle{empty}

%%%%%%%%%%%%%%%%%%%%%%%%%%%%%%%%%%%%%%%%%%%%%%%%%%%%%%%%%%%%%%%%%%%%%%%%%%%%%%%%

\begin{abstract}
% This paper studies the problem of leader-follower continuum deformation coordination of a multi-quadcopter system (MQS) using cooperative localization. While leaders define the desired continuum deformation coordination based on their own positions, followers acquire the desired coordination by following desired trajectories that are defined based on estimations of the leaders' positions. 
% Unmanned aerial vehicles (UAVs) have been widely used in military, non-military{\color{blue},} and academic research applications. Most of UAV systems rely on GPS signals, however outage of GPS signals are inevitable in different scenarios. 
This paper integrates cooperative localization with continuum deformation coordination of a multi-quadcopter system {(MQS) to assure safety and optimality of the quadcopter team coordination in the presence of position uncertainty. We first consider the MQS as a finite number of particles of a deformable triangle in a $3$-D motion space and  define their continuum deformation coordination as a leader-follower problem in which leader quadcopters} can estimate (know) their positions but follower quadcopters  rely on relative position measurements to localize themselves and estimate the leaders' positions. 
We  then propose a navigation strategy for the MQS to plan and acquire the desired continuum deformation coordination, in the presence of measurement noise, disturbance, and position uncertainties, such that collision is avoided and rotor angular speeds of all quadcopters remain bounded. We show the efficacy of the proposed strategy by simulating the continuum deformation coordination of an MQS with eight quadcopters.

\end{abstract}
 
%62467
%%%%%%%%%%%%%%%%%%%%%%%%%%%%%%%%%%%%%%%%%%%%%%%%%%%%%%%%%%%%%%%%%%%%%%%%%%%%%%%%
\section{Introduction}
Unmanned vehicles have been widely used in military~\cite{peng2009design} and non-military applications such as data acquisition from hazardous environments \cite{argrow2005uav} or agricultural farm fields~\cite{tsouros2019review}, traffic surveillance applications~\cite{puri2005survey}, urban search and rescue~\cite{surmann2019integration}, wildlife monitoring and exploration~\cite{witczuk2018exploring} and delivery tasks~\cite{arbanas2016aerial}. % One of the main challenges for a group of UAVs to accomplish the objective is to acquire the accurate global position information. 
Global position estimation is a challenging problem  for unmanned vehicles navigating in uncertain environments. Researchers have proposed feature-based \cite{bailey2006simultaneous} and landmark-based \cite{betke1997mobile} simultaneous localization and mapping (SLAM) algorithms~ for mobile robot localization {in unknown environments}.
% For example, has been widely used for localization in unknown environments. Ref. \cite{betke1997mobile}  uses landmarks to localize mobile robots with respect to the environment.
% Using Global Positioning System (GPS) for navigation is another useful method for localization. Since these methods rely on the existence of distinct features, landmarks or GPS signals, localization cannot be performed in case of outage of GPS signals or lack of distinct features that can be revisited. 
For multi-agent localization, cooperative localization (CL) has been proposed to enable mobile agents to  estimate their global positions by sharing odometry and relative position information. CL has been used in a wide variety of applications such as navigation of double-integrator multi-agents systems \cite{rastgoftar2021continuum} and  ground and aerial vehicles~\cite{sharma2009vision}, search and rescue missions~\cite{jennings1997cooperative}, and target tracking problems~\cite{sharma2012bearing}.

In CL, each agent is equipped with sensors, processing and communication capabilities which enables it to take relative measurements with respect to in-neighbor agents and distribute information to the fusion Center (FC) or only to the in-neighbor agents. These information are mostly noisy signals due to the measurement noises and dynamics of the system. CL uses different estimation approaches, such as extended Kalman filters (EKFs)~\cite{kia2016cooperative}, maximum likelihood~\cite{howard2002localization}, maximum a posteriori (MAP)~\cite{nerurkar2009distributed}, to estimate global positions of member agents of a team by filtering the relative position measurements provided in a distributed fashion.

In this work, we combine CL and continuum deformation coordination approach \cite{rastgoftar2021safe,rastgoftar2018safe}  to safely plan the group coordination of a multi-quadcopter system (MQS) in the presence of position uncertainty.   We consider a group of MQS moving in a 3-D motion space with the desired coordination defined by a non-singular deformation mapping called \textit{homogeneous transformation}. Homogeneous deformation coordination is defined as a leader-follower problem; an $n$-D continuum deformation  of a quadcopters are guided by $n+1$ leader agents, located at vertices of a $n$-D simplex for all time $t$ ($n\in\{1,2,3\}$ denotes the dimension of the continuum deformation coordination).  In this work, without loss of generality, quadcopters are considered as particles of a $2$-D deformable body coordinating in an obstacle-laden motion space, thus, $n=2$, and the desired continuum deformation coordination is defined by three leaders. While the existing homogeneous transformation coordination \cite{rastgoftar2021safe,rastgoftar2018safe} model quadcopters with deterministic dynamics, this paper studies continuum deformation coordination of the MQS in the presence of position uncertainty, measurement noise, and disturbance. In particular, we assume leaders can localize themselves with respect to the environment but followers localize themselves, estimate leaders position, and acquire their desired trajectories by cooperative localization. While the MQS continuum deformation coordination is planned such that travel distance and time are minimized in an obstacle-laden environment, we formally specify and verify safety of the MQS continuum deformation in the presence of global position uncertainty to assure  angular speed of no quadcopter violates a certain upper limit, and collision is avoided.

The organization of the paper is as follows. Section~\ref{se: Problem Statement} presents the problem formulation. Section~\ref{sec: Collective Dynamics of MQS} presents the collective dynamics of MQS. Section~\ref{sec: State Estimation of MQS} presents the state estimation approach and KF. Section \ref{Sec: path planning} discusses the continuum deformation planning in the presence of position uncertainty. %Section VI proposes a path planning and design approach for desired trajectory of the leaders. 
Section~\ref{Sec: Simulation}  gives the simulation of the proposed method on a network of 8 quadcopters, and  followed by Conclusion in Section~\ref{sec: conclusion}.

% papers investigates safety and optimally of continuum deformation coordination of the MQS in the presence position uncertainty, measurement noise, and disturbance. 

% We assign the desired trajectory of the leaders based on the time-optimization problem in which the traveling time between the initial and final configuration is minimized subject to two safety conditions. We suppose that the leaders are equipped with sensors to receive GPS signals, and followers are equipped with proprioceptive sensors to acquire relative pose measurements respect to in-neighbor agents. All agents are connected to an FC to transmit information. 

% Considering the aforementioned setting, we develop a framework for continuum deformation coordination of MQS through simultaneous cooperative localization. We provide the collective dynamics of the quadcopters in which the input is the leaders desired trajectory, and the output only contains the estimated global states of the leaders and the estimated relative states of the followers respect to in-neighbor agents.

%%%%%%%%%%%%%%%%%%%%%%%%%%%%%%%%%%%%%%%%%%%%%%%%%%%%%%%%%%%%%%%%%%%%%%%%%%%%%%%%
\section{Problem Statement}\label{se: Problem Statement}
\begin{figure}
     \centering
     \begin{subfigure}[b]{0.5\textwidth}
         \centering
         \includegraphics[width=\textwidth]{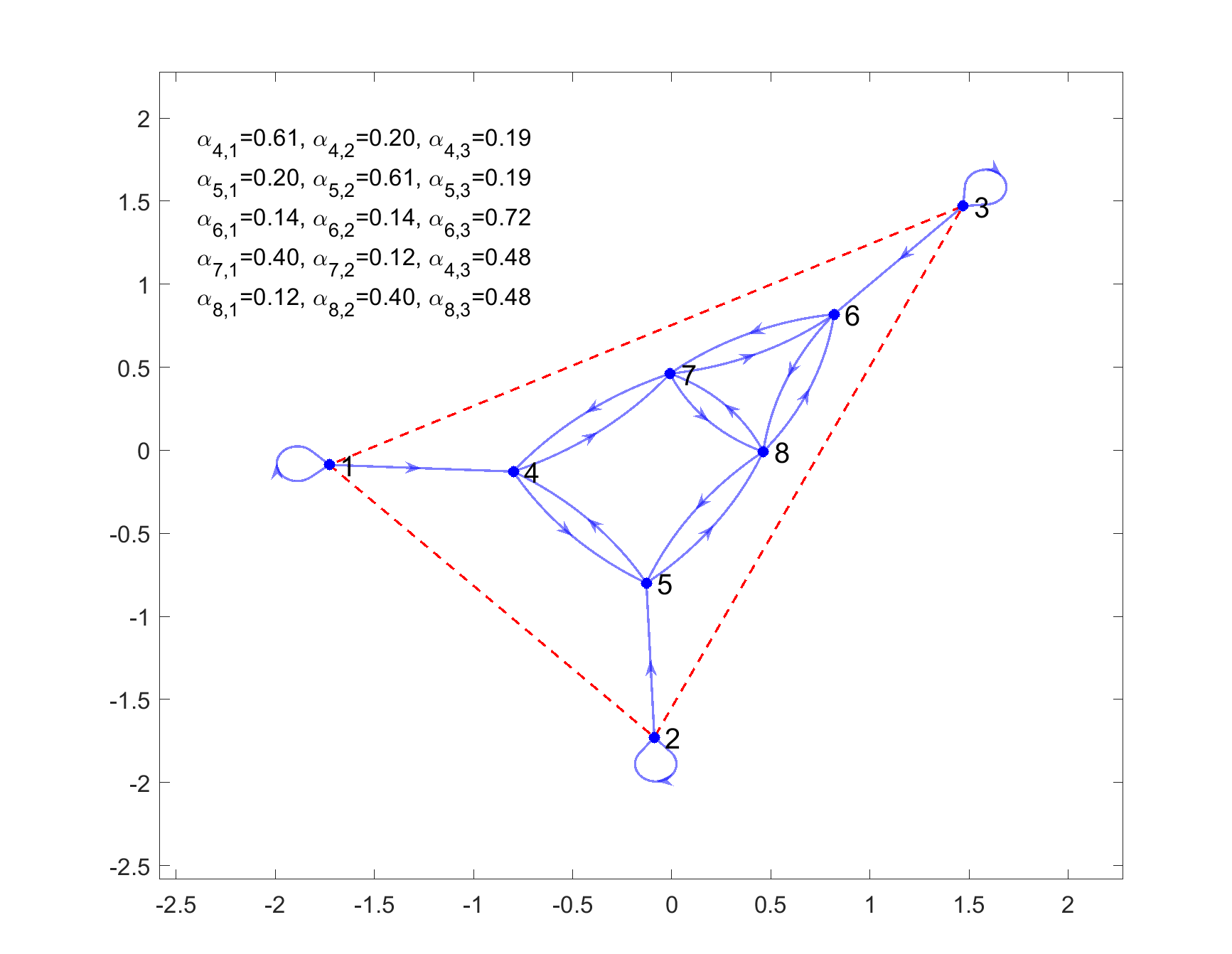}
         \vspace{-0.5cm}
         \caption{}
         \label{fig: graph}
     \end{subfigure}
     \hfill
     \begin{subfigure}[b]{0.45\textwidth}
         \centering
         \includegraphics[width=\textwidth]{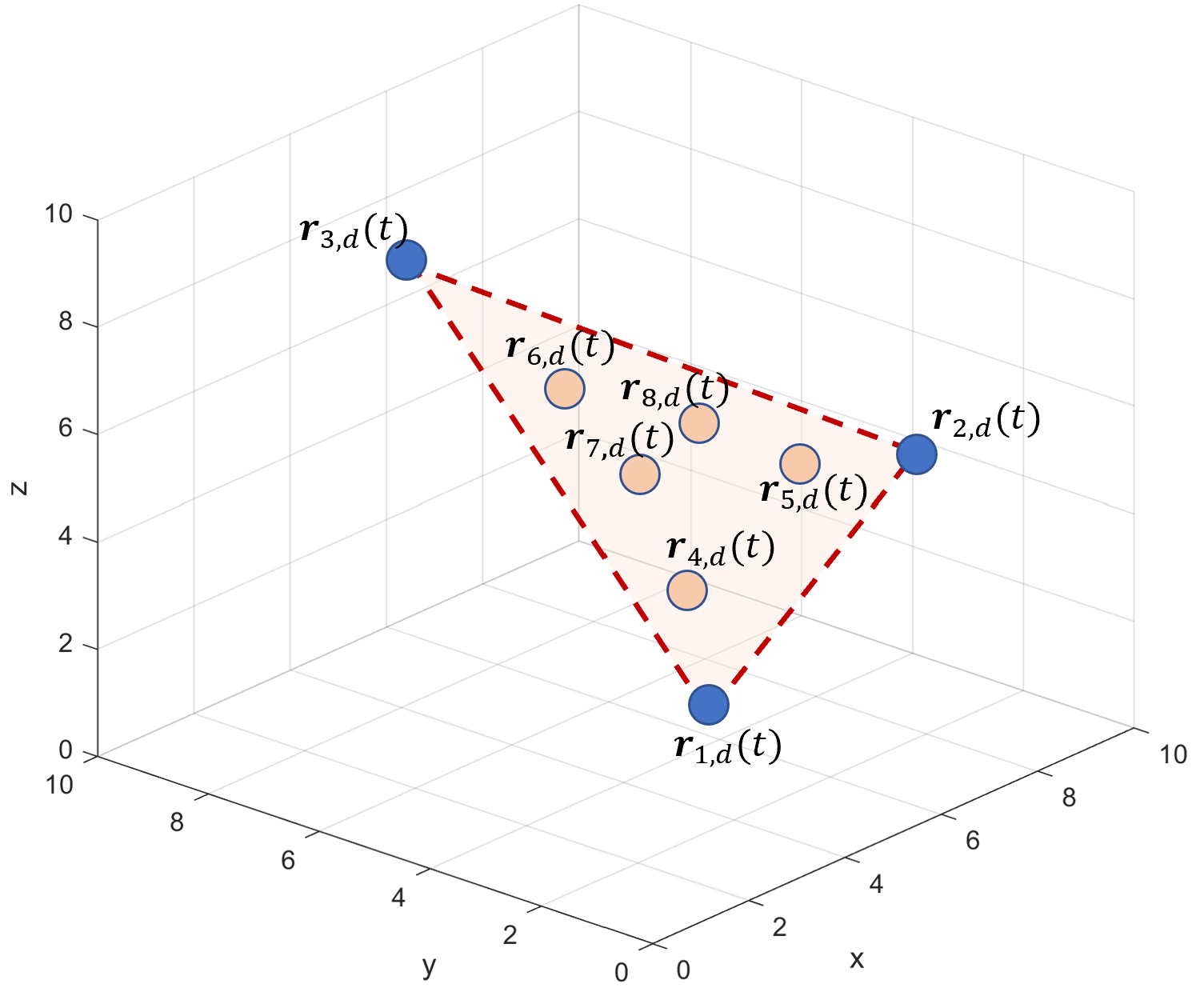}
         \vspace{-0.5cm}
         \caption{}
         \label{fig: desired trajectories}
     \end{subfigure}
     
        \caption{(a) Blue arrows show the directed graph of MQS, and red dashed lines show the leaders' reference configuration in the $x-y$ plane . $\alpha_{i,j}$ are shown in the plot. (b) Agents' configuration on 2-D simplex in 3-D space at time $t$}
        \label{fig: Errors}
\end{figure}

We consider collective motion of an MQS in an obstacle-laden motion space where quadcopters are identified by unique index numbers identified by set $\mathcal{V}=\left\{1,\cdots,N\right\}$. We treat the quadcopters as particles of a $2$-D deformable triangle with three leaders defined by $\mathcal{V}_L=\left\{1,2,3\right\}$ and $N-3$ followers defined by set $\mathcal{V}_F=\mathcal{V}\setminus \mathcal{V}_L=\left\{4,\cdots,N\right\}$. Suppose that an MQS is represented by a directed graph $G(\mathcal{V},\mathcal{E})$ (see Fig~\ref{fig: graph}), where $\mathcal{V}$ is the node set, and the edge set $\mathcal{E} \subseteq \mathcal{V} \times \mathcal{V}$ is defined as a set of pairs $(i,j)$ connecting node $i$ to node $j$ ($i,j\in \mathcal{V}$). Specifically, edge $(i,j)$ physically means that agent $j$ can take the relative measurement of agent $i$. It should be noted that self loop in the network implies that the corresponding agent can receive GPS signals and can measure its global states.  Without loss of generality, we assume that each follower has 3 in-neighbor agent in the network. We denote the in-neighbor agents of agent $i$ as set $\mathcal{N}_i=\{i_1,i_2,i_3\}$.

Let $\mathbf{r}_i(t)=\begin{bmatrix} x_i(t)&y_i(t)&z_i(t) \end{bmatrix}^T$ and $\mathbf{r}_{i,d}(t)=\begin{bmatrix} x_{i,d}(t)&y_{i,d}(t)&z_{i,d}(t) \end{bmatrix}^T$  denote the global position and desired position vector of agent $i\in \mathcal{V}$ at time $t$, respectively. Let $\mathbf{r}_{i,0}=\begin{bmatrix} x_{i,0}&y_{i,0}&0 \end{bmatrix}^T$ denote the reference position of agent $i\in \mathcal{V}$ in $x-y$ plane. 

% For the MQS, set 
% \[
% \Omega_0=\left\{\mathbf{r}_{i,0}=\left(x_{i,0},y_{i,0}\right):\forall i\in \mathcal{V}\right\}
% \]

% defines the reference configuration of the MQS in the $x-y$ plane, where $\mathbf{r}_{i,0}$ is the reference position of quadcopter $i$. 

We assume that the desired trajectory of each quadcopter $i \in \mathcal{V}_L$ is obtained from the following equation:
\begin{eqnarray}\label{equ: desired trajectory}
\mathbf{r}_{i,d}(t) = \mathbf{Q}(t,t_0) \left( \mathbf{r}_{i,0}-\mathbf{d}(t_0) \right)+\mathbf{d}(t)\quad t\in \left[ t_0,t_f\right]
\end{eqnarray}
where $\mathbf{Q}(t,t_0)\in \mathbb{R}^{3\times 3}$ is the Jacobian matrix and $\mathbf{d}(.)\in \mathbb{R}^{3}$ is the rigid body displacement vector~\cite{rastgoftar2018safe}.

We assume that the leaders' desired trajectories are given (see section~\ref{Sec: path planning}), and we define the desired trajectory of the followers as a weighted summation of leaders' position. For every quadcopter $i\in \mathcal{V}_F$, we define three parameters $\alpha_{i,1}$, $\alpha_{i,2}$, and $\alpha_{i,3}$ ($\sum_{j=1}^{3}{\alpha_{i,j}}=1$),  based on reference position of quadcopter $i$ and the leaders' reference positions as follows:
\begin{equation}\label{eq: alpha}
    \begin{bmatrix}
    \alpha_{i,1}\\
    \alpha_{i,2}\\
    \alpha_{i,3}\\
    \end{bmatrix}
    =
    \begin{bmatrix}
    x_{1,0}&x_{2,0}&x_{3,0}\\
    y_{1,0}&y_{2,0}&y_{3,0}\\
    1&1&1\\
    \end{bmatrix}
    ^{-1}
    \begin{bmatrix}
    x_{i,0}\\
    y_{i,0}\\
    1\\
    \end{bmatrix}
    ,\quad \forall i\in \mathcal{V}_F.
\end{equation}
 
The collective motion of the MQS is defined as a leader-follower problem in which the desired  trajectory of quadcopter $i\in \mathcal{V}_F$, denoted by $\mathbf{r}_{i,d}$, are given by
\begin{equation}
    \mathbf{r}_{i,d}(t)=\sum_{j\in\mathcal{V}_L}\alpha_{i,j}\mathbf{r}_{j,d}(t).
\end{equation}

In Fig.~\ref{fig: graph}, $\alpha_{i,j}$'s are represented for a group of 8 quadcopters. Directed graph $G(\mathcal{V},\mathcal{E})$ is also shown in Fig.~\ref{fig: graph}. Fig.~\ref{fig: desired trajectories} shows the desired configuration of an MQS on 2-D simplex in 3-D space at $t$.

\begin{figure*}[h]
\centering
\includegraphics[width=0.9\textwidth]{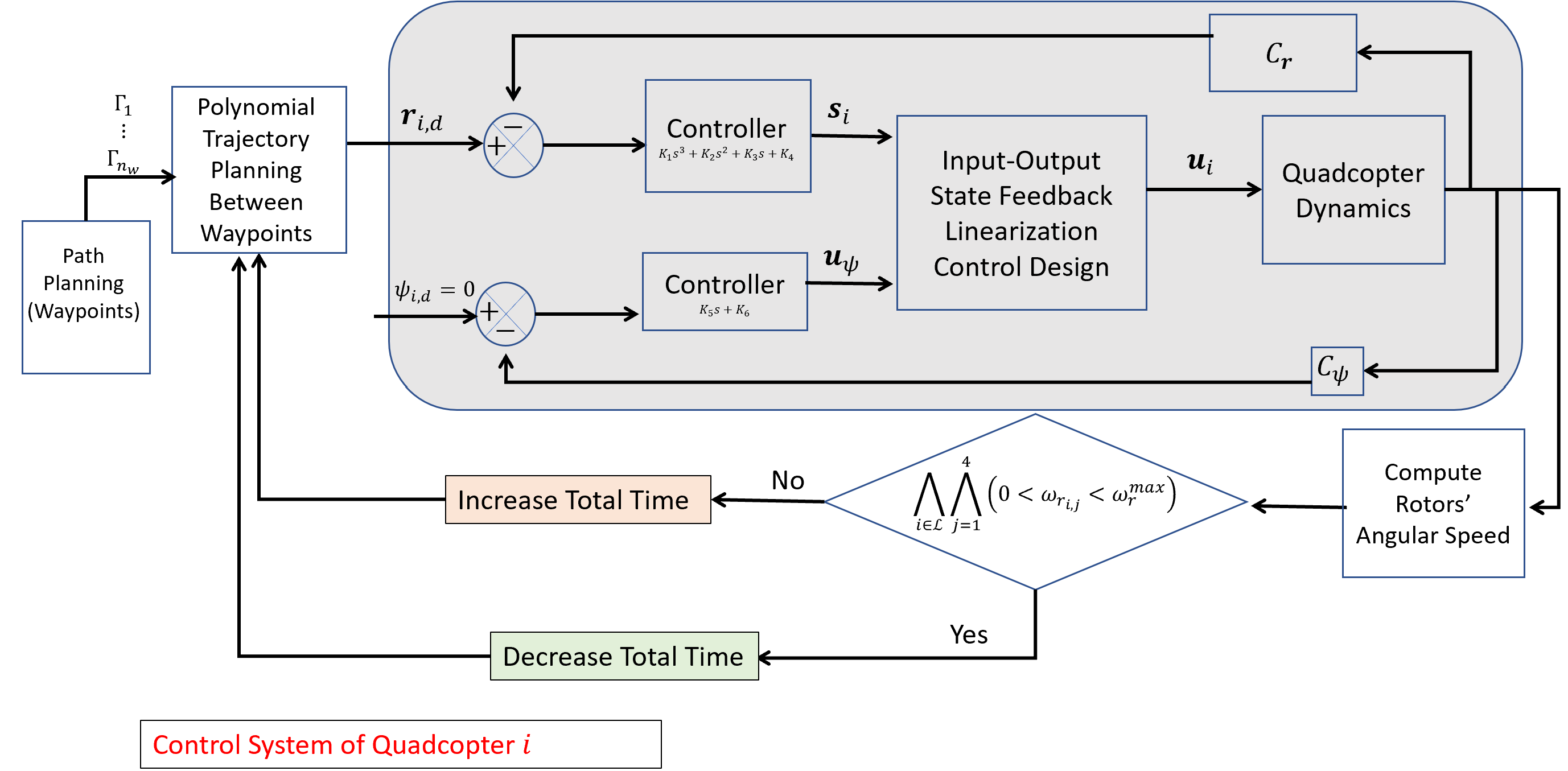}
\caption{Block diagram of MQS coordination control with the proposed method}
\label{Fig: block diagram}
\end{figure*}

The main objective of this work is to design distributed coordination control for an MQS to safe travel in an obstacle-laden environment (see Fig.~\ref{Fig: block diagram}). We suppose that the leaders have access to the GPS signals and followers can only measure the relative pose respect to in-neighbor agents. We consider two safety conditions. First, we assume that the rotor speeds of every quadcopter must not exceed $\omega^{max}_r$. This safety condition can be formally specified by

\begin{equation}\label{safety omega}
    0<{\omega_r}_{i,j}(t)\leq \omega^{max}_r,\quad \forall i\in \mathcal{V},~j\in \left\{1,\cdots,4\right\},~\forall t\geq t_0
\end{equation}
where ${\omega_r}_{i,j}(t)$ is the angular speed of rotor $j\in \left\{1,\cdots,4\right\}$ of quadcopter $i\in \mathcal{V}$ at time $t\geq t_0$. Moreover, trajectory of each quadcopter should be close enough to the corresponding desired trajectory. In particular, following condition should be satisfied for all agents $i\in \mathcal{V}$ at all $t$:
\begin{eqnarray}\label{safety r}
    ||\mathbf{r}_{i}(t)-\mathbf{r}_{i,d}(t)|| < \delta \quad \forall i\in \mathcal{V},~\forall t\geq t_0
\end{eqnarray}
where $\delta$ is the distance threshold value between desired and actual trajectories.  
% where the reference trajectory of leaders are known and determined 

%Given above problem setup, the first objective of this paper is to prescribe the desired trajectory of every leader $j\in \mathcal{V}_L$ such that $\mathbf{r}_{j,d}(t)$ the MQS travel distance is minimized and MQS safety is assured.

\section{Collective Dynamics of MQS}\label{sec: Collective Dynamics of MQS}
In this section, we present the collective dynamics of an MQS in 3-D space. 
 We consider the motion of MQS as particles of a 2-D continuum deformable body guided by 3 leaders. We assume that each agent $i\in \mathcal{V}_L$ is equipped with proprioceptive sensors that can measure the global position. Moreover, we assume that each  agent $i\in \mathcal{V}_F$ can only measure the relative position with respect to the in-neighbor agents.

We assume that the leaders know the desired trajectories of~\eqref{equ: desired trajectory}, and we define the desired trajectory of the followers as a weighted summation of leaders' position in the following form:  
\begin{eqnarray}
\mathbf{r}_{i,d}(t) = \sum_{j=1}^{3}{\alpha_{i,j}\mathbf{r}_{j}(t)}, \quad t\in \left[ t_0,t_f\right] \quad \forall i\in \mathcal{V}_F
\end{eqnarray}
where $\alpha_{i,1},\alpha_{i,2}$ and $\alpha_{i,3}$ are positive numbers associated to agent $i$, and defined in~\eqref{eq: alpha}.

Consequently, weight matrix $\mathbf{W}\in \mathbb{R}^{N\times N}$, defined based on the position of the quadcopters, can be written as 
\begin{eqnarray}
    \mathbf{W} = \left\{ \begin{matrix}
        w_{i,j}&i\in \mathcal{V}_F, j\in \mathcal{V}_L\\
        0&\text{otherwise}
    \end{matrix}\right. .
\end{eqnarray}

From the above definition, matrix $\mathbf{W}$ can be partitioned in the form of
\begin{eqnarray}
    \mathbf{W} = \begin{bmatrix}
    \mathbf{0}_{3\times 3} & \mathbf{0}_{3\times N-3}\\
    \mathbf{W}_0 & \mathbf{0}_{N-3\times N-3}
    \end{bmatrix}
\end{eqnarray}
where $\mathbf{W}_0\in \mathbb{R}^{N-3 \times 3}$ is defined as
\begin{eqnarray}
\mathbf{W}_0 =    \begin{bmatrix}
    \alpha_{4,1} & \alpha_{4,2} & \alpha_{4,3}\\
    \vdots & \vdots & \vdots\\
    \alpha_{N,1} & \alpha_{N,2} & \alpha_{N,3}
    \end{bmatrix}.
\end{eqnarray}

We define matrix $\mathbf{L}$ as 
\begin{eqnarray}
    \mathbf{L} = \mathbf{W}-\mathbf{I}_N.
\end{eqnarray}

Let $\mathbf{X} =\text{vec}\left(\begin{bmatrix}
\mathbf{r}_1 \dots \mathbf{r}_N
\end{bmatrix}^T \right)$ be the concatenation of position vector of all agents. The external dynamics of all quadcopters (see equation~\eqref{external linear dynamics}) can be written in the following form:
\begin{eqnarray}\label{equ: network state space}
    \frac{d}{dt}\left( \begin{bmatrix}
    \mathbf{X}\\
    \Dot{\mathbf{X}}\\
    \ddot{\mathbf{X}}\\
    \dddot{\mathbf{X}}
    \end{bmatrix}\right)
    =\mathbf{A}_{SYS}\begin{bmatrix}
    \mathbf{X}\\
    \Dot{\mathbf{X}}\\
    \ddot{\mathbf{X}}\\
    \dddot{\mathbf{X}}
    \end{bmatrix}
    + \mathbf{B}_{SYS}\begin{bmatrix}
    \mathbf{S}_L\\
    \Dot{\mathbf{S}}_L\\
    \ddot{\mathbf{S}}_L\\
    \dddot{\mathbf{S}}_L
    \end{bmatrix}\\ \label{equ: network state space 2}
    \mathbf{Y} = \mathbf{C}_{SYS} \begin{bmatrix}
    \mathbf{X}^T&\Dot{\mathbf{X}}^T&\ddot{\mathbf{X}}^T&\dddot{\mathbf{X}}^T
    \end{bmatrix}^T
\end{eqnarray}
where $\mathbf{A}_{SYS}\in \mathbb{R}^{12N\times 12N},\mathbf{B}_{SYS}\in \mathbb{R}^{12N\times 12N}$  and $\mathbf{C}_{SYS}\in \mathbb{R}^{36(N-2)\times 12N}$ are defined as
\begin{eqnarray}
    \mathbf{A}_{SYS} = \begin{bmatrix}
    \mathbf{0}_{3N\times 3N} & \mathbf{I}_{3N} & \mathbf{0}_{3N\times 3N} & \mathbf{0}_{3N\times 3N} \\
    \mathbf{0}_{3N\times 3N} & \mathbf{0}_{3N\times 3N} & \mathbf{I}_{3N} & \mathbf{0}_{3N\times 3N}\\
    \mathbf{0}_{3N\times 3N} & \mathbf{0}_{3N\times 3N} & \mathbf{0}_{3N\times 3N} & \mathbf{I}_{3N}\\
    K_4 \mathbf{I}_3\otimes \mathbf{L} & K_3 \mathbf{I}_3\otimes \mathbf{L} & K_2 \mathbf{I}_3\otimes \mathbf{L} & K_1 \mathbf{I}_3\otimes \mathbf{L}
    \end{bmatrix}
\end{eqnarray}
\begin{eqnarray}
    \mathbf{B}_{SYS} = \begin{bmatrix}
    \mathbf{0}_{9N\times 9} & \mathbf{0}_{9N\times 9} & \mathbf{0}_{9N\times 9} & \mathbf{0}_{9N\times 9}\\
    K_4 \mathbf{I}_3\otimes \mathbf{L}_0 & K_3 \mathbf{I}_3\otimes \mathbf{L}_0 & K_2 \mathbf{I}_3\otimes \mathbf{L}_0 & K_1 \mathbf{I}_3\otimes \mathbf{L}_0
    \end{bmatrix}
\end{eqnarray}
\begin{eqnarray}
    \mathbf{C}_{SYS} = \mathbf{I}_{12} \otimes \mathbf{C}_{0}
\end{eqnarray}
where $\mathbf{L}_0\in \mathbb{R}^{N\times 3}$ is
\begin{eqnarray}
    \mathbf{L}_0 = \begin{bmatrix}
    \mathbf{I}_3\\
    \mathbf{0}_{N-3\times 3}
    \end{bmatrix}.
\end{eqnarray}
$\mathbf{S}_L$ is defined as concatenation of desires trajectories of leaders as follows:
\begin{eqnarray}
    \mathbf{S}_L= \text{vec}\left(\begin{bmatrix}
    \mathbf{r}_{1,d} &\mathbf{r}_{2,d}& \mathbf{r}_{3,d}
    \end{bmatrix}^T \right)
\end{eqnarray}

$\mathbf{C}_0\in \mathbb{R}^{3(N-2)\times N}$ is a matrix with the $ij^\text{th}$ entry ${C_0}_{i,j}$ defined in the following way:
\begin{eqnarray}
    {C_0}_{i,j} =  \left\{ \begin{matrix}
        1&i=j , j\in \mathcal{V}_L\\
        1& (j,i)\in \mathcal{E}, j\in \mathcal{V}_F\\
        -1& i=j, j\in \mathcal{V}_F\\
        0&\text{otherwise}
    \end{matrix}\right. .
\end{eqnarray}

We define $\mathbf{Y}_d = \text{vec} \left( \begin{bmatrix}
\mathbf{r}_{1,d} & \dots & \mathbf{r}_{N,d}
\end{bmatrix}^T\right)$. Vector $\mathbf{Y}_d$ and $\mathbf{S}_L$ are related as 
\begin{eqnarray}
    \mathbf{Y}_d = (\mathbf{I}_3 \otimes \mathbf{H})\mathbf{S}_L.
\end{eqnarray}
where $\mathbf{H} = -\mathbf{L}^{-1}\mathbf{L}_0$~\cite{rastgoftar2021safe} .

Now, by defining $\mathbf{E}(t)= \mathbf{Y}(t)-\mathbf{Y}_d(t)$, the error dynamics can be written in the form of
\begin{eqnarray}\label{equ: network Error}
    \frac{d}{dt}\left( \begin{bmatrix}
    \mathbf{E}\\
    \Dot{\mathbf{E}}\\
    \ddot{\mathbf{E}}\\
    \dddot{\mathbf{E}}
    \end{bmatrix}\right)
    =\mathbf{A}_{SYS}\begin{bmatrix}
    \mathbf{E}\\
    \Dot{\mathbf{E}}\\
    \ddot{\mathbf{E}}\\
    \dddot{\mathbf{E}}
    \end{bmatrix}
    + \begin{bmatrix}
    \mathbf{0}\\
    \mathbf{0}\\
    \mathbf{0}\\
    \mathbf{I}_3 \otimes \mathbf{H}^T
    \end{bmatrix}\ddddot{\mathbf{S}}_L
\end{eqnarray}

%%%%%%%%%%%%%%%%%%%%%%%%%%%%%%%%%%%%%%%%%%%%%%%%%%%%%%%%%%%%%%%%%%%%%%%%%%%%%%%%
\section{State Estimation of MQS}\label{sec: State Estimation of MQS}
In this section, we present the Kalman Filter state estimation algorithm following \cite{simon2006optimal}. We consider a centralized scenario in which all agents share their measurements to a FC. Collective dynamics of the system is presented in~\eqref{equ: network state space}. Note that leaders can measure their states, and followers can only measure the relative states of the in-neighbor agents. Matrix $\mathbf{C}_0$ in~\eqref{equ: network state space 2} represents the explicit form of the absolute and the relative measurements in the network $G$. 

In the first step, we discretize the continuous dynamics of~\eqref{equ: network state space}. Suppose that sensors are sampling every $\Delta t$ second.  Discretizing the continuous state space model~\eqref{equ: network state space} and ~\eqref{equ: network state space 2} lead to the following discrete approximation model 

\begin{eqnarray}\label{equ: discrete network state space}
\mathbf{x}_{[k+1]} &=& (\mathbf{A}_{SYS}\Delta t +\mathbf{I}) \mathbf{x}_{[k]} + (\mathbf{B}_{SYS}\Delta t)\mathbf{u}_{[k]}+\mathbf{\eta}_{[k]}\\
    \mathbf{y}_{[k+1]} &=& \mathbf{C}_{SYS} \mathbf{x}_{[k]}+\mathbf{\nu}_{[k]}
\end{eqnarray}
where ${\mathbf{x}}_{[k]}$, ${\mathbf{u}}_{[k]}$ and ${\mathbf{y}}_{[k]}$ represent the state vector, the control input vector and the measurement vector at time-step $k$, respectively, in in~\eqref{equ: network state space},\eqref{equ: network state space 2}.  $\mathbf{\eta}_{[k]}$ and $\mathbf{\nu}_{[k]}$ are process noise and measurement noise, respectively. We assume that $\mathbf{\eta}_{[k]}$ and $\mathbf{\nu}_{[k]}$ are zero-mean independent white Gaussian processes with known covariances $\mathbf{Q}_{[k]}$, $\mathbf{R}_{[k]}$, respectively.  For each time step the Kalman filter is given by the following expressions:
\begin{eqnarray}
\mathbf{P}^-_{[k+1]} &=& (\mathbf{A}_{SYS}\Delta t +\mathbf{I}) \mathbf{P}^+_{[k]} (\mathbf{A}\Delta t +\mathbf{I})^T + \mathbf{Q}_{[]k]}\\
\mathbf{K}_{[k+1]} &=& \mathbf{P}^-_{[k+1]} \mathbf{C}_{SYS}^T\left(\mathbf{C}_{SYS}\mathbf{P}^-_{[k+1]}\mathbf{C}_{SYS}^T+\mathbf{R}_{[k+1]}\right)^{-1}\\
\mathbf{x}^-_{[k+1]} &=& (\mathbf{A}_{SYS}\Delta t +\mathbf{I}) \mathbf{x}^+_{[k]} + (\mathbf{B}_{SYS}\Delta t)\mathbf{u}_{[k]}\\
\mathbf{x}^+_{[k+1]} &=& \mathbf{x}^-_{[k+1]} + \mathbf{K}_{[k+1]}\left(\mathbf{y}_{[k+1]}-\mathbf{C}\mathbf{x}^-_{[k+1]} \right)\\
\mathbf{P}^+_{[k+1]} &=& \left( \mathbf{I}-\mathbf{K}_{[k+1]}\mathbf{C}_{SYS}\right) \mathbf{P}^-_{[k+1]}
\end{eqnarray}
where $``^+"$,$``^-"$ refer to the prior and posterior estimation, respectively. That is to say, $``^+",``^-"$ correspond to the estimation after and before we process the measurement at time step $k$, respectively. $\mathbf{P}_{[k]},\mathbf{K}_{[k]}$ represent the error estimation covariance and Kalman filter gain at time step $k$, respectively.
%%%%%%%%%%%%%%%%%%%%%%%%%%%%%%%%%%%%%%%%%%%%%%%%%%%%%%%%%%%
\section{Path Planning}\label{Sec: path planning}
\begin{figure}[ht]
\centering
\includegraphics[width=0.5\textwidth]{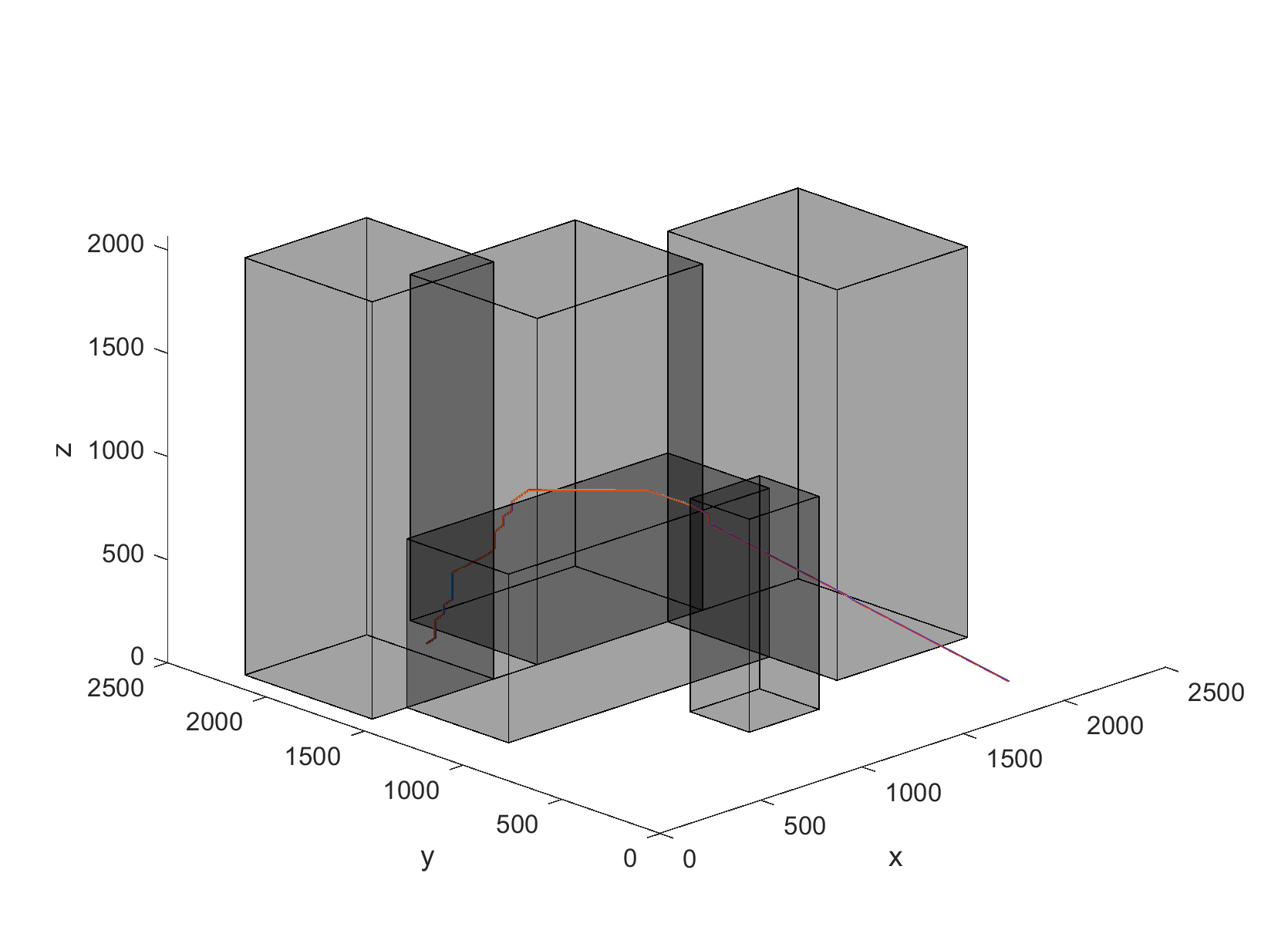}
\caption{An obstacle-laden environment. Leaders' desired paths which is generated from the approach discussed in Section \ref{Sec: path planning}, also shown in the plot }
\label{Fig: obstacles}
\end{figure}

We use the A* search method~\cite{hart1968formal} to path planning of the leaders in an obstacle-laden environment (see  Fig.\ref{Fig: obstacles}). Deploying A* search method results to a line-graph in which the node set represents the waypoints, and edge set represents the path segment between waypoints. We denote the position of waypoints by $\Gamma_1,\dots,\Gamma_n$. We assume that each quadcopter starts with zero velocity and zero acceleration at start point $\Gamma_i$, and reaches to the end point $\Gamma_{i+1}$ of each segment with zero velocity and zero acceleration. We consider a trajectory of a quadcopter as a polynomial function of time with zero velocity and acceleration at $\Gamma_i$ and $\Gamma_{i+1}$. This results to the desired trajectories of the leaders which satisfy~\eqref{equ: desired trajectory} as follows:
\begin{eqnarray}
    \mathbf{r}^j_{i,d}(t) = (1-\beta(t))\Gamma_j + \beta(t)\Gamma_{j+1}
\end{eqnarray}
where $\beta(t) = \frac{6}{T^5_j}t^5 - \frac{15}{T^4_j}t^4 + \frac{10}{T^3_j}t^3$, and superscript $j$ in $\mathbf{r}^j_{i,d}(t)$ denotes the $j^\text{th}$ path segment between $\Gamma_j$ and $\Gamma_{j+1}$. We denote the total travelling time by $T$; we linearly allocate travelling time $T_j$ to the path segment between waypoints $\Gamma_j$ and $\Gamma_{j+1}$ based on the travelling distance between $\Gamma_j$ and $\Gamma_{j+1}$. From~\eqref{equ: network Error}, the tracking error can be written as
\begin{eqnarray}\label{equ: network Error 2}
    \begin{bmatrix}
    \mathbf{E}\\
    \Dot{\mathbf{E}}\\
    \ddot{\mathbf{E}}\\
    \dddot{\mathbf{E}}
    \end{bmatrix}
    =\int_{t_0}^{t} e^{\mathbf{A}_{SYS}(t-\eta)}\begin{bmatrix}
    \mathbf{0}\\
    \mathbf{0}\\
    \mathbf{0}\\
    \mathbf{I}_3 \otimes \mathbf{H}^T
    \end{bmatrix}\ddddot{\mathbf{S}}_L \,d\eta 
\end{eqnarray}
From the above expression, as $T_j$ tends to infinity, $\mathbf{E}$ also tends to 0. This leads to the fact that there exists an optimal time $T^*$ for give $\delta$ such that safety condition~\eqref{safety r} is satisfied for all time.

In order to find the optimal traveling time for MQS subjected to that safety conditions~\eqref{safety omega} and~\eqref{safety r}, we use bisection method. We initiate with a large $T$ such that~\eqref{safety omega} and~\eqref{safety r} are satisfied. Using the bisection method, we keep updating $T$ until one of the safety conditions is violated. We denote the optimal time by $T^*$. 
%%%%%%%%%%%%%%%%%%%%%%%%%%%%%%%%%%%%%%%%%%%%%%%%%%%%%%%%%%%
\section{Simulation}\label{Sec: Simulation}
In this section, we consider an MQS containing 8 quadcopters labeled as $\mathcal{V}=\{1,\dots,8\}$. In order to acquire the continuum deformation coordination, We consider 3 leaders in this group, labeled as $\mathcal{V}_L=\{1,2,3\}$, and the rest of agents are considered as followers $\mathcal{V}_F=\{4,\dots,8\}$. A directed graph $G(\mathcal{V},\mathcal{E})$ is generated based on proximity for local relative measurements (see Fig.~\ref{fig: graph}). We assume that leaders are equipped by proprioceptive sensors which enable to acquire the global state vector measurements at each time step. On the other hand, each follower can only measure the relative sates with respect to its neighbors (e.g. agent 4 can take the measurements relative to agents 1, 5 and 7). Note that self-loop in the network implies that the corresponding agent can measure its global states. Quadcopters' specification are listed in Table~\ref{Table}.

\begin{table}[h!]
\centering
 \begin{tabular}{|c |c| c| c|} 
 \hline
 $m$ & $g$ & $l$ & $I_x$ \\  
 \hline
 0.468 & 9.81 & 0.225 & $4.856\times 10^{-3}$\\
 \hline \hline
 $I_y$ & $I_z$ & $b$ & $k$\\ 
 \hline
 $4.856\times 10^{-3}$ & $8.801\times 10^{-3}$ & $2.98\times 10^{-6}$ & $1.14\times 10^{-7}$\\
 \hline
 \end{tabular}
 \caption{Quadcopters' specification}
 \label{Table}
\end{table}

\begin{figure}[h]
\centering
\includegraphics[width=0.4\textwidth]{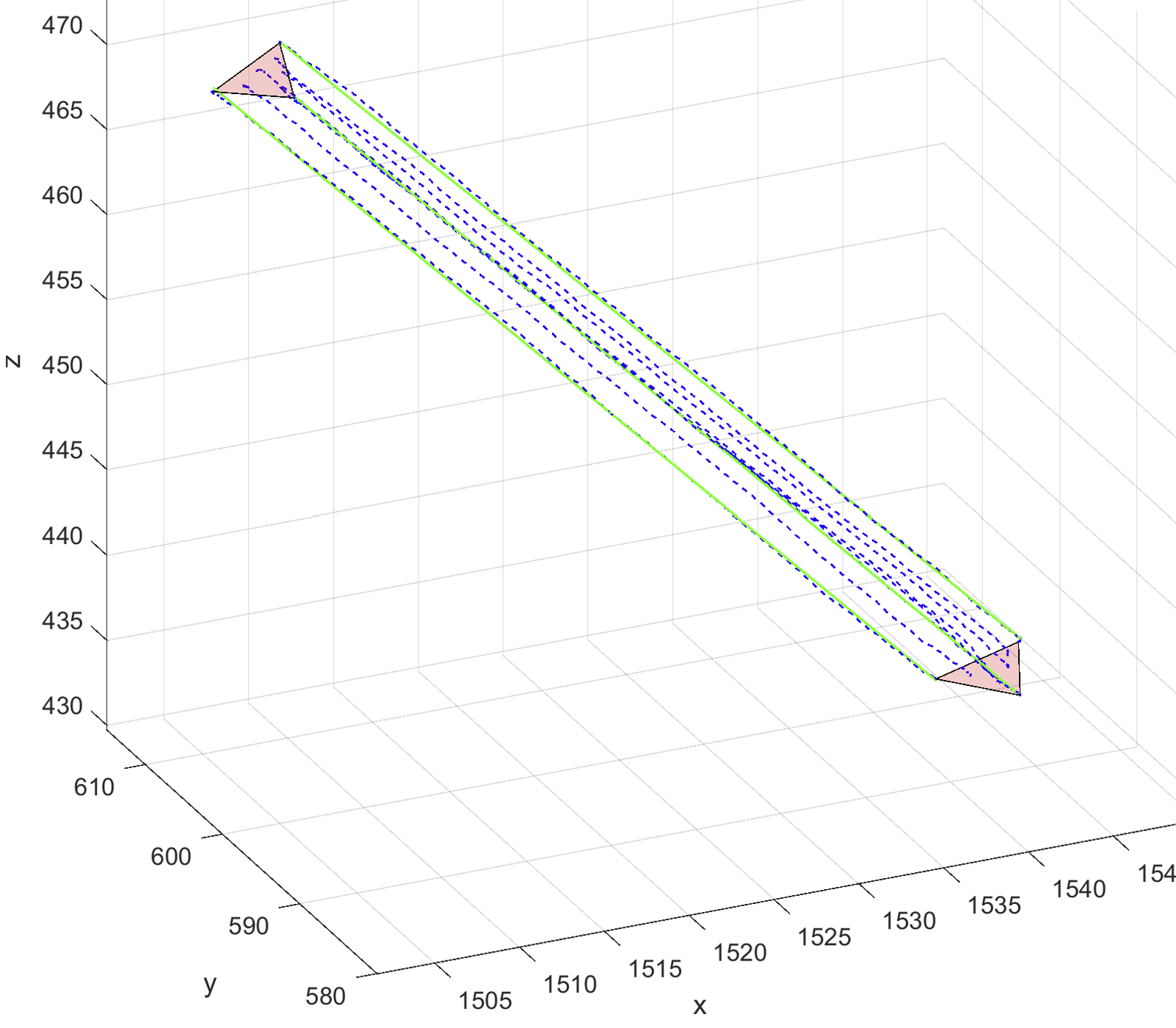}
\caption{Actual (green lines) and desired (blue lines) trajectories of MQS for the path segment between waypoints $\Gamma_2,\Gamma_3$.}
\label{Fig: trajectories_1}
\end{figure}

We consider the standard deviation of 0.1 for process noise $Q$ and measurement noise $R$. Sampling time in our simulation is 0.01 $\sec$. Fig.~\ref{Fig: trajectories_1} shows the trajectories of the MQS from $\Gamma_2$ to $\Gamma_3$. We choose $K_1=10,K_2=35,K_3=50$ and $K_4=35$. Blue dashed lines show the actual trajectories of the agents, and solid green lines show the desired trajectories of 3 leaders. As shown in Fig.~\ref{Fig: trajectories_1}, followers are contained in the triangle, formed by the three leaders. Fig.~\ref{fig: Errors} shows the estimation and tracking error of agent 4 (a follower agent).   

\begin{figure}[h]
\centering
\includegraphics[width=0.5\textwidth]{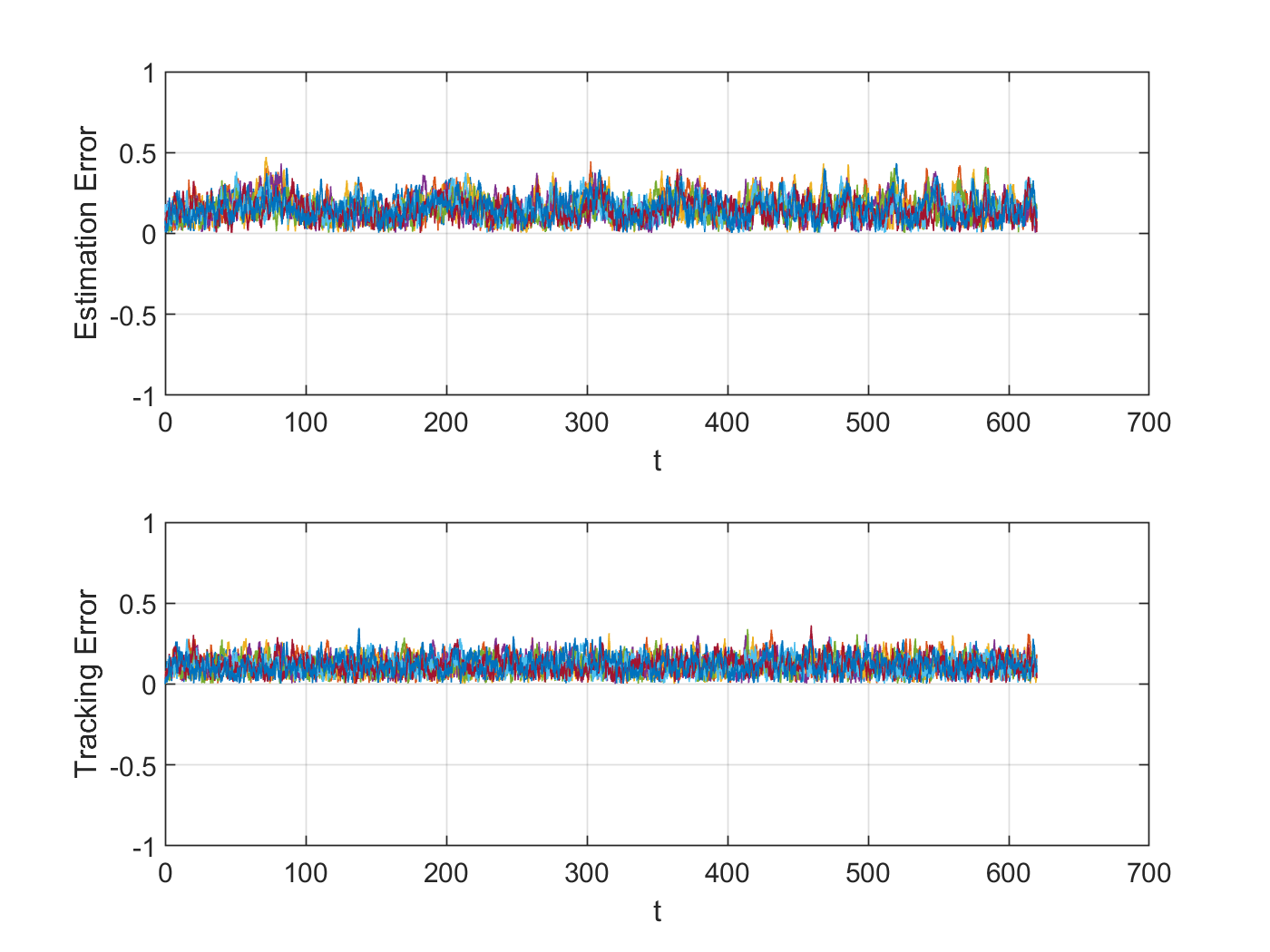}
\caption{Tracking error and estimation error for all agents}
\label{Fig: error}
\end{figure}

\begin{figure}[h]
\centering
\includegraphics[width=0.4\textwidth]{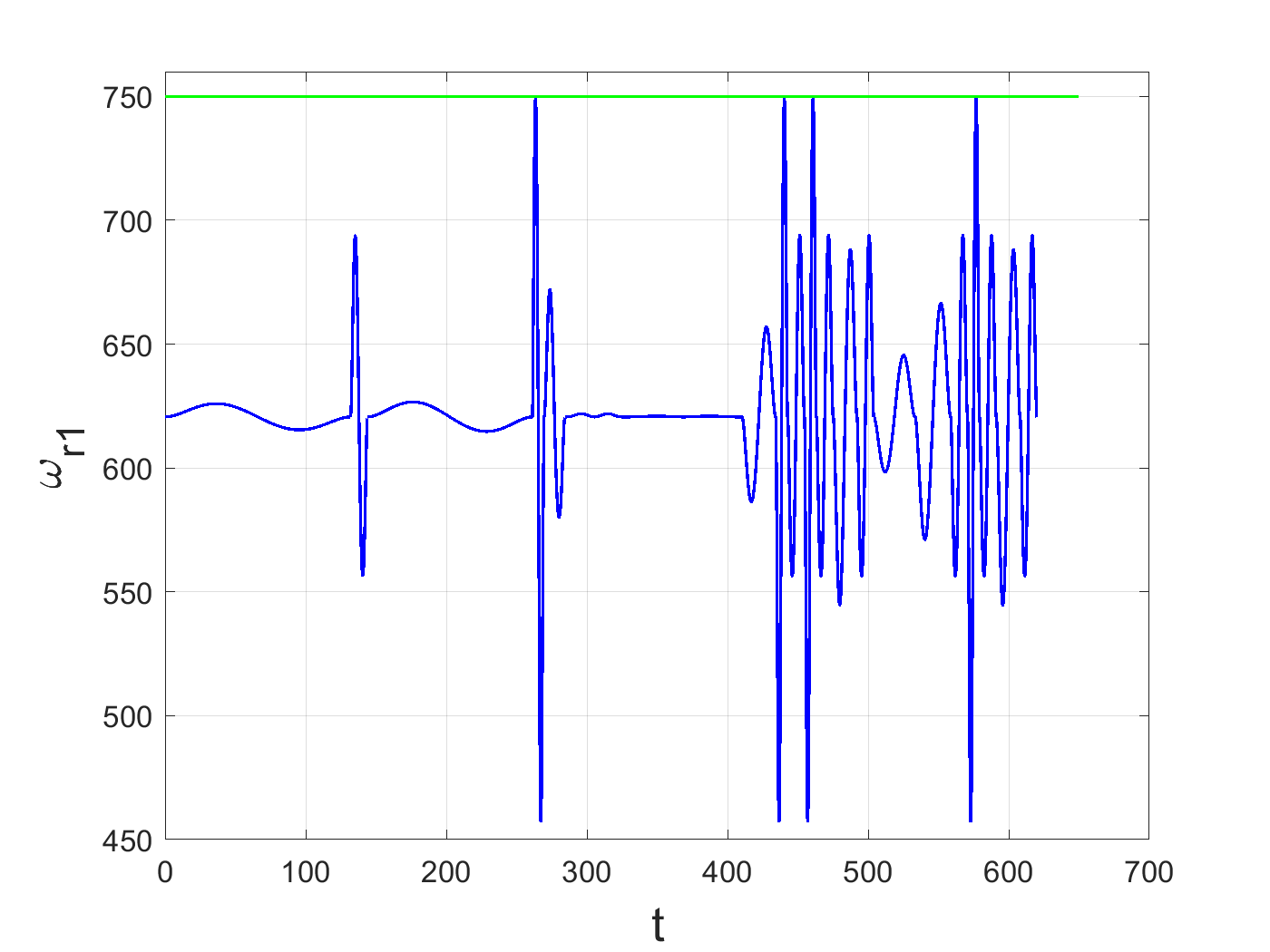}
\caption{Blue line shows the angular speed ${\omega_r}_1$ during the total travelling time. Green line shows the upper limit ${\omega_r}^{\max}$.}
\label{Fig: angular_speed}
\end{figure}

\begin{figure}[h]
\centering
\includegraphics[width=0.4\textwidth]{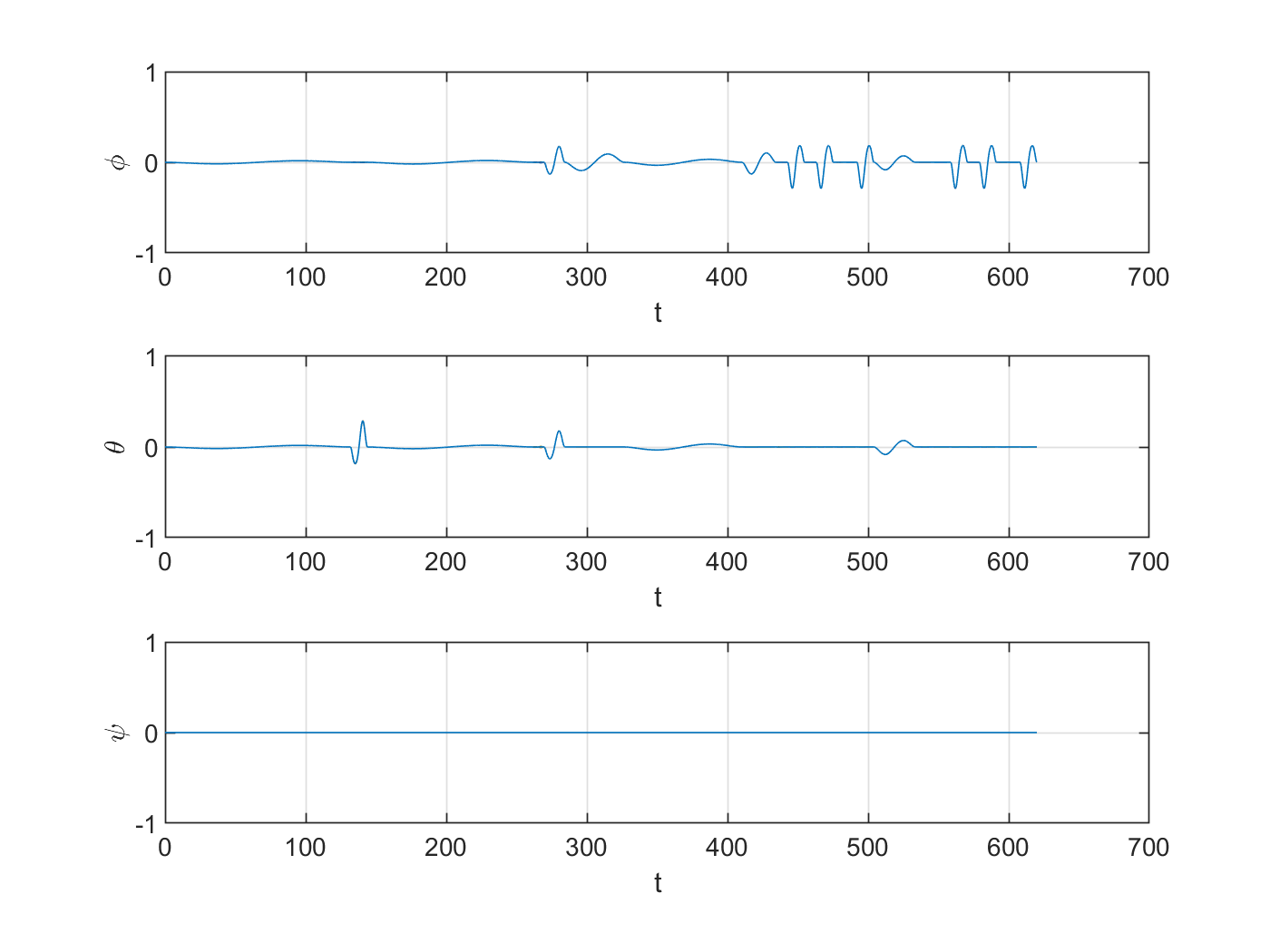}
\caption{$\phi , \theta$ and $\psi$ for agent 4.}
\label{Fig: attitude}
\end{figure}

\begin{figure}[h]
\centering
\includegraphics[width=0.4\textwidth]{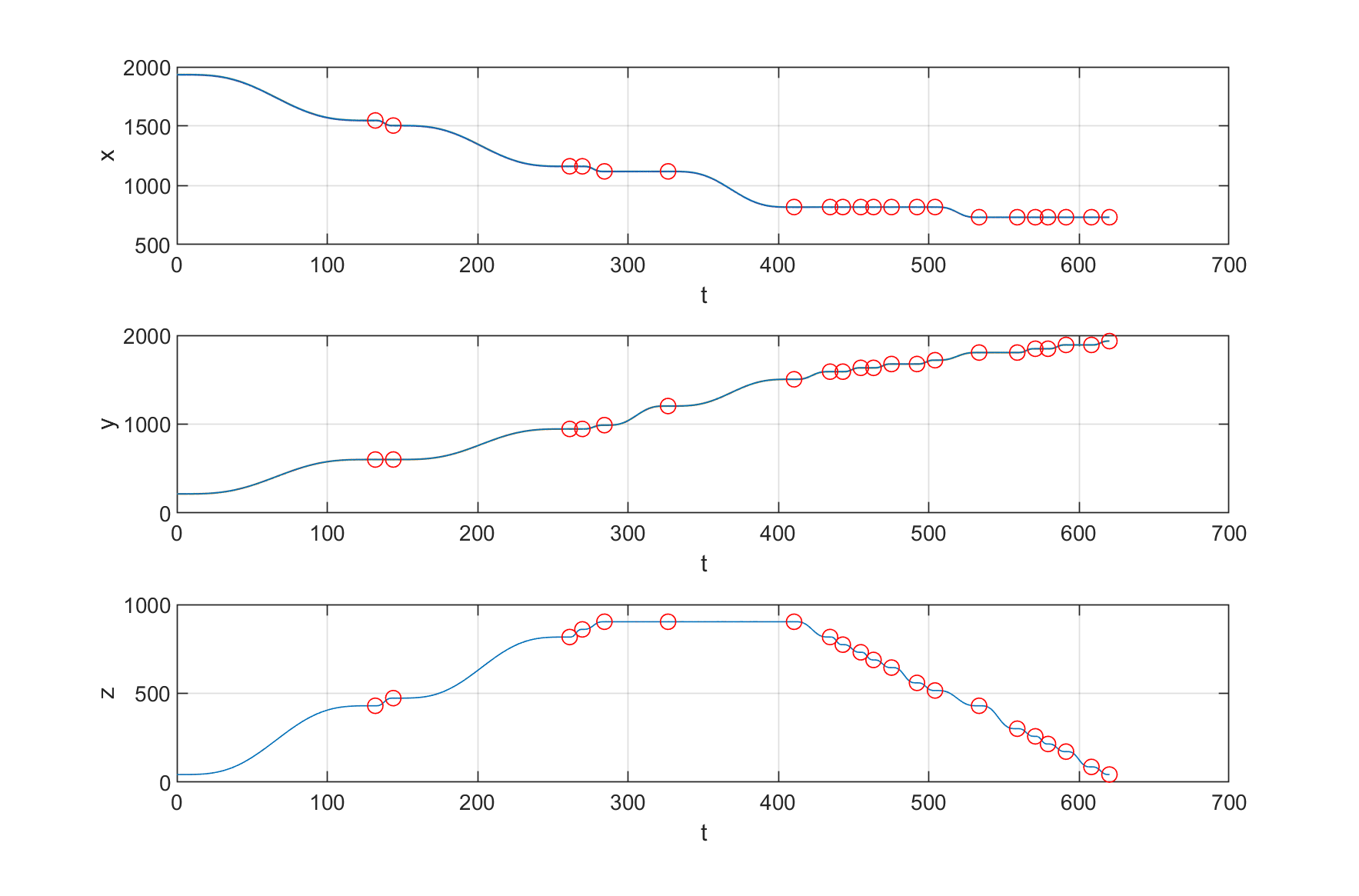}
\caption{$x,y$ and $z$ components of quadcopters. Red circles show the waypoints.}
\label{Fig: xyz}
\end{figure}

 We choose $\delta =0.5$ in the safety condition~\eqref{safety r}, and ${\omega_r}_{\max}=750$ in safety condition~\eqref{safety omega}. Using the bisection method, $T^*$ is computed as 640 seconds. Fig.~\ref{Fig: angular_speed} shows the angular speed of rotor 1  for quadcopter 4. As shown in Fig.\ref{Fig: angular_speed}, $\omega_r$ is not exceeding the ${\omega_r}_{\max}=750$. Fig.~\ref{Fig: attitude} and~\ref{Fig: xyz} show the roll, pitch and yaw angle of quadcopter 4 and $x,y$ and $z$ components of quadcopters for $t\in [0,T^*]$, respectively.
%%%%%%%%%%%%%%%%%%%%%%%%%%%%%%%%%%%%%%%%%%%%%%%%%%%%%%%%%%%
\section{Conclusion}\label{sec: conclusion}
We  developed  a framework for continuum deformation coordination of MQS through  simultaneous  cooperative  localization.  We  provided the collective dynamics of the quadcopters in which the input is the leader's desired trajectory, and the output only contains the  estimated  global  states  of  the  leaders  and  the  estimated relative states of the followers respect to in-neighbor agents. We used Kalman Filter for state estimation of the collective motion system. In this work, we used FC to collect and distribute the information to the network. As a part of the future work we plan to develop a decentralized method for state estimation and coordination of quadcopters. 
%%%%%%%%%%%%%%%%%%%%%%%%%%%%%%%%%%%%%%%%%%%%%%%%%%%%%%%%%%
\section{Acknowledgement}
This work has been supported by the National Science
Foundation under Award Nos. 2133690 and 1914581.

%%%%%%%%%%%%%%%%%%%%%%%%%%%%%%%%%%%%%%%%%%%%%%%%%%%%%%%%%%%
\appendix

\subsection{Mathematical Modeling of Quadcopters and Trajectory Tracking Control}
\label{Mathematical Modeling of Quadcopters and Trajectory Tracking Control}

In this work, we consider the following assumptions in mathematical modeling of quadcopter motions.

\begin{assumption}
Quadcopter is a symmetrical rigid body with respect to the axes of body-fixed frame.
\end{assumption}
\begin{assumption}
Aerodynamic loads are neglected due to low speed assumption for quadcopters.
\end{assumption}

Let $\hat{\boldsymbol{i}},\hat{\boldsymbol{j}},\hat{\boldsymbol{k}}$ be the base unit vectors of inertial coordinate system, and $\hat{\boldsymbol{i}}_b,\hat{\boldsymbol{j}}_b,\hat{\boldsymbol{k}}_b$ be the base unit vectors of a body-fixed coordinate system whose origin is at the center of mass of the quadcopter. In this section, for convenience, we omitted $i$ subscript of $i^{\text{th}}$ quadcopter in the governing equations. The attitude of the quadcoper is defined by three Euler angles $\phi,\theta$ and $\psi$ as roll angle, pitch angle and yaw angle, respectively. In this work, we use 3-2-1 standard Euler angles to determine orientation of the quadcopter. Therefore, the rotation matrix between fixed-body frame and the inertial frame can be written as
\begin{eqnarray}
    \boldsymbol{R}(\phi,\theta,\psi) = \boldsymbol{R}(\phi,0,0)\boldsymbol{R}(0,\theta,0)\boldsymbol{R}(0,0,\psi)\\
    =\begin{bmatrix}
     c_{\theta}c_{\psi}&  s_{\theta}c_{\psi}s_{\phi}-s_{\psi}c_{\phi}& s_{\theta}c_{\psi}c_{\phi}+s_{\psi}s_{\phi}\\ 
     c_{\theta}s_{\psi}&  s_{\theta}s_{\psi}s_{\phi}+c_{\psi}c_{\phi}& s_{\theta}s_{\psi}c_{\phi}-c_{\psi}s_{\phi}\\ 
    -s_{\theta}&  c_{\theta}s_{\phi}& c_{\theta}c_{\phi}
    \end{bmatrix}
\end{eqnarray}
where $s(.)=\sin(.),c(.)=\cos(.)$. Let $\boldsymbol{r}=\begin{bmatrix} x&y&z \end{bmatrix}^T$ denote the position of the center of mass of the quadcopter in inertial frame, and $\boldsymbol{\omega}=\begin{bmatrix}
\omega_x&\omega_y&\omega_z
\end{bmatrix}^T$  denote the angular velocity of the quadcopter represented in the fixed-body frame.

Using the Newton-Euler formulas, equations of motion of a quadcopter can be written in the following form:
\begin{eqnarray}\label{translation_eq}
    &\Ddot{\boldsymbol{r}} = -g\hat{\boldsymbol{k}}+\frac{p}{m}\hat{\boldsymbol{k}}_b,\\\label{rotation_eq}
    &\Dot{\boldsymbol{\omega}} = -\boldsymbol{J}^{-1}\left[\boldsymbol{\omega}\times(\boldsymbol{J}\boldsymbol{\omega})\right]+\boldsymbol{J}^{-1}\boldsymbol{\tau},
\end{eqnarray}
where $m,\boldsymbol{J}$ denote, respectively, mass and mass moment of inertia of the quadcopter. $g$ is the gravity acceleration and $p$ is the thrust force generated by the four rotors. Relation between the thrust force $p$ and angular speed of the rotors, denoted by $\omega_{r_i}$, can be written as
\begin{eqnarray}
    p = \sum_{i=1}^{4}{{f_r}_i} = b\sum_{i=1}^{4}{{{\omega^2_{r_i}}}},
\end{eqnarray}
where $b$ is the aerodynamic force constant ($b$ is a function of the density of air, the shape of the blades, the number of the blades, the chord length of the blades, the pitch angle of the blade airfoil and the drag constant), and ${f_r}_i$ is the thrust force of $i^{\text{th}}$ rotor. In~\eqref{translation_eq}, $\boldsymbol{\tau}=\begin{bmatrix} \tau_\phi&\tau_\theta&\tau_\psi
\end{bmatrix}^T$ is the control torques generated by four rotors. Relation between the $\boldsymbol{\tau}$ and angular speed of the rotors can be written in the following form
\begin{eqnarray}
    \boldsymbol{\tau}=\begin{bmatrix}
     \tau_\phi \\ 
    \tau_\theta \\ 
    \tau_\psi
    \end{bmatrix}=
    \begin{bmatrix}
     bl({\omega^2_{r_4}}-{\omega^2_{r_2}}) \\ 
    bl({\omega^2_{r_3}}-{\omega^2_{r_1}}) \\ 
    k({\omega^2_{r_2}}+{\omega^2_{r_4}}-{\omega^2_{r_1}}-{\omega^2_{r_3}})
    \end{bmatrix},    
\end{eqnarray}
where $l$ is the distance of each rotor from center of the quadcopter, and $k$ is a positive constant corresponding to the aerodynamic torques. By concatenating $p$ and $\boldsymbol{\tau}$ as input vector to the system, we can write
\begin{eqnarray} \label{eq: angular rotor speed}
    \boldsymbol{u}=\begin{bmatrix}
    p\\
     \tau_\phi \\ 
    \tau_\theta \\ 
    \tau_\psi
    \end{bmatrix}=\begin{bmatrix}
    b & b & b & b\\
    0 & -bl & 0 & bl\\ 
    -bl & 0 & bl & 0 \\ 
    -k & k & -k & k
    \end{bmatrix}\begin{bmatrix}
    {\omega^2_{r_1}}\\
    {\omega^2_{r_2}} \\ 
    {\omega^2_{r_3}} \\ 
    {\omega^2_{r_4}}
    \end{bmatrix}.
\end{eqnarray}

By defining state vector $\boldsymbol{x}=\begin{bmatrix}\boldsymbol{r}^T&\Dot{\boldsymbol{r}}^T&\phi&\theta&\psi&\boldsymbol{\omega}^T
\end{bmatrix}^T$ and input vector $\boldsymbol{u}=\begin{bmatrix}
p& \tau_{\phi}&\tau_{\theta}&\tau_{\psi}
\end{bmatrix}^T$, \eqref{translation_eq},\eqref{rotation_eq} can be written in the state space non-linear form of 
\begin{eqnarray} \label{non-linear state space}
     \left\{\begin{matrix}
        \Dot{\boldsymbol{x}}= \boldsymbol{f}(\boldsymbol{x})+\boldsymbol{g}(\boldsymbol{x})\boldsymbol{u}\\ 
        \boldsymbol{y}=\boldsymbol{C}\boldsymbol{x}
    \end{matrix}\right.   
\end{eqnarray}
where, $\boldsymbol{f}(\boldsymbol{x})$ and $\boldsymbol{g}(\boldsymbol{x})$ are defined as 
\begin{eqnarray}
     \boldsymbol{f}(\boldsymbol{x})=\begin{bmatrix}
    \boldsymbol{v}\\
     -g\hat{\boldsymbol{k}} \\ 
    \boldsymbol{\Gamma^{-1}}\boldsymbol{\omega} \\ 
    -\boldsymbol{J}^{-1}\left[\boldsymbol{\omega}\times(\boldsymbol{J}\boldsymbol{\omega})\right]
    \end{bmatrix},
\end{eqnarray}

\begin{eqnarray}
    \boldsymbol{g}(\boldsymbol{x})=\begin{bmatrix}
    \boldsymbol{0}_{3\times1} & \boldsymbol{0}_{3\times3}\\
    \frac{\hat{\boldsymbol{k}}_b}{m} & \boldsymbol{0}_{3\times3}\\ 
    \boldsymbol{0}_{3\times1} & \boldsymbol{0}_{3\times3}\\ 
    \boldsymbol{0}_{3\times1} & \boldsymbol{J}^{-1}
    \end{bmatrix}
\end{eqnarray}
and $\boldsymbol{C} = [\boldsymbol{I}_{3\times3}, \boldsymbol{0}_{3\times9}]$. $\boldsymbol{v}$ is the velocity vector of the quadcopter, and $\boldsymbol{\Gamma}$ is the matrix which relates Euler angular velocity to the angular velocity of the quadcopter. $\boldsymbol{0}_{i\times j}$ is a $i\times j $ zero matrix. In order to find $\boldsymbol{\Gamma}$, we can represent $\boldsymbol{\omega}$ in the following form
\begin{eqnarray}\label{angular velocity}
    \boldsymbol{\omega} = \Dot{\psi} \hat{\boldsymbol{k}}_1+\Dot{\theta} \hat{\boldsymbol{j}}_2 + \Dot{\phi} \hat{\boldsymbol{i}}_b,
\end{eqnarray}
where $\hat{\boldsymbol{j}}_2=\boldsymbol{R}(\phi,0,0)\hat{\boldsymbol{j}}_b$ and $\hat{\boldsymbol{k}}_1=\boldsymbol{R}(\phi,\theta,0)\hat{\boldsymbol{k}}_b$. Consequently, 
\begin{eqnarray}
    \boldsymbol{\Gamma} =\begin{bmatrix}
    1 & 0 &-s_{\theta}\\
    0 & c_{\phi} &c_{\theta}s_{\phi}\\ 
    0 & -s_{\phi} &c_{\phi}c_{\theta}
    \end{bmatrix} .
\end{eqnarray}
From~\eqref{angular velocity}, the angular acceleration $\Dot{\boldsymbol{\omega}}$ can be formulated in the following way:
\begin{eqnarray}\label{angular velocity 1}
    \Dot{\boldsymbol{\omega}} = \Tilde{\boldsymbol{B}}_1 \begin{bmatrix}
    \ddot{\phi}& \ddot{\theta}&\ddot{\psi}
    \end{bmatrix}^T +\Tilde{\boldsymbol{B}}_2
\end{eqnarray}
where $\Tilde{\boldsymbol{B}}_1=\begin{bmatrix}
\hat{\boldsymbol{i}}_b& \hat{\boldsymbol{j}}_2& \hat{\boldsymbol{k}}_1
\end{bmatrix}$ and
\begin{eqnarray} \label{tilde_B2}
    \Tilde{\boldsymbol{B}}_2=\Dot{\theta}\Dot{\psi}(\hat{\boldsymbol{k}}_1 \times \hat{\boldsymbol{j}}_2)+ \Dot{\phi}(\Dot{\psi}\hat{\boldsymbol{k}}_1+\Dot{\theta}\hat{\boldsymbol{j}}_2 )\times \hat{\boldsymbol{i}}_b
\end{eqnarray}
On the other hand, from~\eqref{non-linear state space},
\begin{eqnarray} \label{angular velocity 2}
    \Dot{\boldsymbol{\omega}} = \boldsymbol{J}^{-1}\left( -\boldsymbol{\omega}\times(\boldsymbol{J}\boldsymbol{\omega}) + \begin{bmatrix}
    u_{2}& u_3 &u_4
    \end{bmatrix}^T \right).
\end{eqnarray}
From~\eqref{angular velocity 1} and \eqref{angular velocity 2}
\begin{eqnarray}\label{equation: u to ddot}
    \begin{bmatrix}
    u_2\\ u_3\\ u_4
    \end{bmatrix} = \boldsymbol{J}\Tilde{\boldsymbol{B}}_1\begin{bmatrix}
     \ddot{\phi}\\ \ddot{\theta}\\\ddot{\psi} 
    \end{bmatrix} + \boldsymbol{J}\Tilde{\boldsymbol{B}}_2+\boldsymbol{\omega}\times(\boldsymbol{J}\boldsymbol{\omega})
\end{eqnarray}

\subsection{Input-Output Feedback Linearization Control}
In this subsection, we provide the input control for the non-linear state space system~\eqref{non-linear state space} to track the desired trajectory $\boldsymbol{r_d}$. We suppose $\boldsymbol{r}_d$ is a smooth function for all  $t\geq t_0$ (i.e. $\boldsymbol{r}_d$ has derivatives of all orders with respect to time). 

In this work, we use the input-output feedback linearization approach\cite{slotine1991applied} to design the input control for a quadcopter to track the desired trajectory~\cite{rastgoftar2021safe}. We use the Lie derivative notation which is defined in the following.
\begin{definition}
Let $h:\mathbb{R}^n\rightarrow  \mathbb{R}$ be a smooth scalar function, and $\boldsymbol{f}:\mathbb{R}^n\rightarrow  \mathbb{R}^n$ be a smooth vector field on $\mathbb{R}^n$. Lie derivative of $h$  with respect to $\boldsymbol{f}$ is a scalar function defined by $L_{\boldsymbol{f}} h=\nabla h\boldsymbol{f}$.
\end{definition}

Concept of input-output linearization is based on differentiating the output until the input appears in the derivative expression. Since $u_2,u_3$ and $u_4$ do not appear in the derivative of outputs, we use the technique, called dynamic extension, in which we redefine the input vector $\boldsymbol{u}$ as the derivative of some of the original system inputs. In particular, we define $\Tilde{\boldsymbol{x}}=\begin{bmatrix}\boldsymbol{x}^T&p&\Dot{p}
\end{bmatrix}^T$ and $\Tilde{\boldsymbol{u}}=\begin{bmatrix}
u_p&\tau_{\phi}&\tau_{\theta}&\tau_{\psi}
\end{bmatrix}^T$. Therefore, extended dynamics of the quadcopter can be expressed in the following form \cite{rastgoftar2021safe}:
\begin{eqnarray} \label{extended non-linear state space}
     \left\{\begin{matrix}
        \Dot{\Tilde{\boldsymbol{x}}}= \Tilde{\boldsymbol{f}}(\Tilde{\boldsymbol{x}})+\Tilde{\boldsymbol{g}}(\Tilde{\boldsymbol{x}})\Tilde{\boldsymbol{u}}\\ 
        \boldsymbol{r}=\Tilde{\boldsymbol{C}}\Tilde{\boldsymbol{x}}
    \end{matrix}\right.   
\end{eqnarray}
where, $\Tilde{\boldsymbol{f}}(\Tilde{\boldsymbol{x}})$ and $\Tilde{\boldsymbol{g}}(\Tilde{\boldsymbol{x}})$ are defined as 
\begin{eqnarray}
     \Tilde{\boldsymbol{f}}(\Tilde{\boldsymbol{x}})=\begin{bmatrix}
    \boldsymbol{f}(\boldsymbol{x})\\
     \Dot{p} \\ 
    0\end{bmatrix}+
    \begin{bmatrix}
    \boldsymbol{0}_{3\times1}\\
     \frac{p}{m}\hat{\boldsymbol{k}}_b \\ 
    \boldsymbol{0}_{8\times1}\end{bmatrix},
\end{eqnarray}

\begin{eqnarray}
    \Tilde{\boldsymbol{g}}(\Tilde{\boldsymbol{x}})=\begin{bmatrix}
    \boldsymbol{0}_{9\times1} & \boldsymbol{0}_{9\times3}\\
    \boldsymbol{0}_{3\times1} & \boldsymbol{J}^{-1}\\ 
    0 & \boldsymbol{0}_{1\times3}\\ 
    1 & \boldsymbol{0}_{1\times3}
    \end{bmatrix}.
\end{eqnarray}
Let $\Tilde{\boldsymbol{g}}_i(\Tilde{\boldsymbol{x}})$ denote the $i^\text{th}$ column of matrix $\Tilde{\boldsymbol{g}}(\Tilde{\boldsymbol{x}})$ and $\Tilde{\boldsymbol{u}}=\begin{bmatrix}\Tilde{u}_1\dots\Tilde{u}_4
\end{bmatrix}^T$ where $\Tilde{u}_1,\dots,\Tilde{u}_4$ corresponds to $u_p,\tau_{\phi},\tau_{\theta},\tau_{\psi}$, respectively. We consider the position of the quadcopter as the output of the system (i.e. $x,y,z$). Inputs appear in the fourth order derivative of the outputs. particularly, for $q\in \{x,y,z\}$ 
\begin{eqnarray}
    \ddddot{q} = L^4_{\Tilde{\boldsymbol{f}}}q + \sum_{i=1}^{4}{L_{\Tilde{\boldsymbol{g}_i}}{L_{\Tilde{\boldsymbol{f}}}}^3q}\Tilde{u}_i
\end{eqnarray}
where $L_{\Tilde{\boldsymbol{g}_i}}{L_{\Tilde{\boldsymbol{f}}}}^3q \neq 0$ for $i=1,\dots,4$. By choosing the state transformation $\mathcal{T}(\Tilde{\boldsymbol{x}})=\begin{bmatrix}
\boldsymbol{z}&\boldsymbol{\zeta}
\end{bmatrix}^T$, \eqref{extended non-linear state space} can be converted to the following internal and external dynamics:
\begin{eqnarray}\label{internal linear dynamics}
    \Dot{\boldsymbol{\zeta}} = \begin{bmatrix}
    0 & 0\\
    0 & 1 \end{bmatrix} \boldsymbol{\zeta} +
    \begin{bmatrix}
    0\\
    1 \end{bmatrix} u_{\psi}
\end{eqnarray}
\begin{eqnarray}\label{external linear dynamics}
    \Dot{\boldsymbol{z}} = \boldsymbol{A}\boldsymbol{z}+\boldsymbol{B}\boldsymbol{s}
\end{eqnarray}
where $\boldsymbol{z} = \begin{bmatrix}
\boldsymbol{r}^T& \Dot{\boldsymbol{r}}^T& \Ddot{\boldsymbol{r}}^T& \dddot{\boldsymbol{r}}^T, 
\end{bmatrix}^T$, and  $\boldsymbol{\zeta} = \begin{bmatrix}
\psi&\Dot{\psi}
\end{bmatrix}^T$
\begin{eqnarray}
    \boldsymbol{A}=\begin{bmatrix}
    \boldsymbol{0}_{9\times 3} & \boldsymbol{I}_{9}\\
    \boldsymbol{0}_{3\times 3} & \boldsymbol{0}_{3\times 9}\end{bmatrix}, \boldsymbol{B}=\begin{bmatrix}
    \boldsymbol{0}_{9\times 3} \\
    \boldsymbol{I}_{3} \end{bmatrix}
\end{eqnarray}
where $\boldsymbol{I}_{j}$ is a $j\times j$ identity matrix. 

Next, we can figure out the Control inputs $\boldsymbol{s}$ and $u_\psi$, such that the linear systems \eqref{internal linear dynamics} and \eqref{external linear dynamics} track the desired trajectory $\boldsymbol{r}_d$. By choosing
\begin{eqnarray}
    u_{\psi} = -K_1\Dot{\psi}-K_2\psi
\end{eqnarray}
where $K_1>0,K_2>0$. Thus, the internal dynamics~\eqref{internal linear dynamics} asymptotically converges to $\boldsymbol{0}$. Moreover, we choose 
\begin{eqnarray}
    \boldsymbol{s}=K_3\left( \dddot{\boldsymbol{r}}_d-\dddot{\boldsymbol{r}}\right)+ K_4\left( \ddot{\boldsymbol{r}}_d-\ddot{\boldsymbol{r}}\right)+ \\ \nonumber
    K_5\left( \dot{\boldsymbol{r}}_d-\dot{\boldsymbol{r}}\right) + K_6\left( {\boldsymbol{r}}_d-{\boldsymbol{r}}\right)
\end{eqnarray}
where $K_3,\dots,K_6$ can be chosen such that the roots of the characteristic equation 
\begin{eqnarray}
    \lambda^4+\lambda^3K_3+\lambda^2K_4+\lambda K_5+K_6 =0,
\end{eqnarray}
are located in the open left half complex plane. Hence, $\boldsymbol{r}$ converges to $\boldsymbol{r}_d$.

In order to find the relation between $\boldsymbol{s}$ and $\Tilde{\boldsymbol{u}}$, we need to find $\ddddot{\boldsymbol{r}}$ by differentiating twice with respect to time from $\ddot{\boldsymbol{r}}$. From~\eqref{non-linear state space}, we have
\begin{eqnarray}
    \ddot{\boldsymbol{r}} = \frac{p}{m}\hat{\boldsymbol{k}}_b-g\hat{\boldsymbol{k}}
\end{eqnarray}
By differentiating the above expression,
\begin{eqnarray}
    \dddot{r} = \frac{\Dot{p}}{m}\hat{\boldsymbol{k}}_b + \frac{p}{m} \boldsymbol{\omega}\times\hat{\boldsymbol{k}}_b
\end{eqnarray}
\begin{eqnarray}\label{ddddotr}
    \ddddot{r} = \frac{1}{m}(\boldsymbol{O}_1\boldsymbol{\Theta}+\boldsymbol{O}_2),
\end{eqnarray}
where $\boldsymbol{\Theta}=\begin{bmatrix}\ddot{p}&\ddot{\phi}&\ddot{\theta}&\ddot{\psi}\end{bmatrix}^T$ and
\begin{eqnarray}
    \boldsymbol{O}_1 = \begin{bmatrix}
     \hat{\boldsymbol{k}}_b& -p\hat{\boldsymbol{j}}_b& p(\hat{\boldsymbol{j}}_2 \times \hat{\boldsymbol{k}}_b)& p(\hat{\boldsymbol{k}}_1)\times \hat{\boldsymbol{k}}_b 
    \end{bmatrix}
\end{eqnarray}
\begin{eqnarray}
    \boldsymbol{O}_2 = p\Tilde{\boldsymbol{B}_2}\times \hat{\boldsymbol{k}}_b + \boldsymbol{\omega}\times(\boldsymbol{\omega\times p \hat{\boldsymbol{k}}_b})+2\Dot{p}\boldsymbol{\omega\times\hat{\boldsymbol{k}}_b}
\end{eqnarray}
where $\Tilde{\boldsymbol{B}_2}$ is defined in~\eqref{tilde_B2}. From~\eqref{equation: u to ddot}, $\boldsymbol{\Theta}$ can be written in the form of
\begin{eqnarray}\label{Theta}
    \boldsymbol{\Theta} = \boldsymbol{O}_3\Tilde{\boldsymbol{u}} + \boldsymbol{O}_4,
\end{eqnarray}
where
\begin{eqnarray}
    \boldsymbol{O}_3 =  
    \begin{bmatrix}
    1 & \boldsymbol{0}_{1\times3}\\
    \boldsymbol{0}_{1\times3} & \boldsymbol{J}^{-1}{\Tilde{\boldsymbol{B}}_1}^{-1} \end{bmatrix},
\end{eqnarray}
\begin{eqnarray}
    \boldsymbol{O}_4 =  
    \begin{bmatrix}
    0\\
    -{\Tilde{\boldsymbol{B}}_1}^{-1}\Tilde{\boldsymbol{B}}_2-\boldsymbol{J}^{-1}\boldsymbol{\omega}\times(\boldsymbol{J}\boldsymbol{\omega}) \end{bmatrix}.
\end{eqnarray}
Substituting \eqref{Theta} in \eqref{ddddotr}
\begin{eqnarray}\label{control input}
    \boldsymbol{s} = \frac{1}{m} \left( \boldsymbol{O}_1\boldsymbol{O}_3\Tilde{\boldsymbol{u}} + \boldsymbol{O}_1\boldsymbol{O}_4+\boldsymbol{O}_2\right)
\end{eqnarray}

\bibliographystyle{IEEEtran}
\bibliography{ref}
\end{document}